\documentclass[12pt,a4paper]{article}

\usepackage[T1]{fontenc}
\usepackage[utf8]{inputenc}
\usepackage[british]{babel}
\usepackage{amsmath,amssymb,amsthm}
\usepackage{url}

\DeclareMathOperator{\Hol}{Hol}

\newtheorem{theorem}{Theorem}[section]

\newtheorem{corollary}[theorem]{Corollary}
\theoremstyle{remark}

\newcommand\numberbracesold{1\,515\,429}
\newcommand\numberbraces{10\,326\,821}
\newcommand\totalbraces{15\,095\,601}
\title{Enumeration of left braces with additive group $C_2\times C_2\times C_4\times C_4$}
\author{A. Ballester-Bolinches\thanks{Departament de Matem\`atiques, Universitat de Val\`encia, Dr.\ Moliner, 50, 46100 Burjassot, Val\`encia, Spain; \texttt{Adolfo.Ballester@uv.es}, \texttt{Ramon.Esteban@uv.es}, \texttt{Vicent.Perez-Calabuig@uv.es}; ORCID 0000-0002-2051-9075, 0000-0002-2321-8139, 0000-0003-4101-8656}\and R. Esteban-Romero\addtocounter{footnote}{-1}\footnotemark\and V. P{\'e}rez-Calabuig\addtocounter{footnote}{-1}\footnotemark}
\begin{document}

\maketitle

\begin{abstract}
  We show that the number of isomorphism classes of left braces of order~$64$ with additive group isomorphic to $C_2\times C_2\times C_4\times C_4$ is $\numberbraces$. This completes the classification of left braces of order~$64$, that turn out to fall into $\totalbraces$ isomorphism classes.

  \emph{Mathematics Subject Classification (2020):}
  16T25, 
  81R50, 
  20-08. 

  \emph{Keywords:} Skew left brace, left brace, regular subgroup, holomorph.
\end{abstract}

\section{Introduction}

A \emph{skew left brace} is a mathematical structure $(B, {+}, {\cdot})$ composed of a set $B$ and two binary operations ${+}$, and ${\cdot}$, on~$B$ for which $(B, {+})$ and $(B, {\cdot})$ are groups linked by the relation $x(y+z)=xy-x+xz$ for $x$, $y$, $z\in B$. Skew left braces were introduced by Guarnieri and Vendramin in \cite{GuarnieriVendramin17}. They generalise \emph{left braces} introduced by Rump in his seminal paper \cite{Rump07}, corresponding to skew left braces with $(B, {+})$ abelian. Here we use $ab$ to denote $a\cdot b$, and  $-a$ to denote the symmetric of $a$ in $(B, {+})$; the expression $a-b$ will denote $a+(-b)$. We also follow the usual hierarchy between the operations $+$ and~$\cdot$.

The determination of the isomorphism classes of (skew) left braces of a given finite order is one of the most natural problems in the study of skew left braces. Guarnieri and Vendramin presented in \cite[Algorithm~5.1]{GuarnieriVendramin17} an algorithm to enumerate the isomorphism classes of skew left braces with a given finite additive group~$A$. The numbers of isomorphism classes of left braces of all orders up to $120$ except $32$, $64$, $81$, and $96$ are presented in \cite{GuarnieriVendramin17}. These numbers were obtained with an implementation in \textsc{Magma} \cite{BosmaCannonFiekerSteel16-Magma} of that algorithm. The computation of the remaining cases, that is, left braces of orders~$32$, $64$, $81$, and $96$, is presented as an open problem in \cite[Problem~6.1]{GuarnieriVendramin17}. Vendramin also left as an open problem (see \cite[Problem~2.13]{Vendramin19-agta}) the construction of all isomorphism classes of left braces of order~$32$. He presented some partial results in \cite[Table~2.3]{Vendramin19-agta}. He considered the numbers of isomorphism classes of left braces of size~$64$, $96$ or $128$ ``extremely large'' and their ``computational methods  not strong enough to construct them all.''

Bardakov, Neshchadim, and Yadav modified in \cite[Algorithm~2.4]{BardakovNeshchadimYadav20}  the method presented in \cite[Algorithm~5.1]{GuarnieriVendramin17} to compute all isomorphism classes of finite skew left braces with a prescribed additive group. They used it to give the number of isomorphism classes of skew left braces of orders~$32$ and~$81$, as well as the number of isomorphism classes of left braces of order~$96$. They were also able to enumerate the isomorphism classes of left braces of order~$64$ in \cite[Table~6]{BardakovNeshchadimYadav20} for all isomorphism classes of abelian additive groups of order~$64$ except for the cases of additive group isomorphic to $C_4\times C_4\times C_4$ (\texttt{SmallGroup(64, 55)} in the notation of the library of small groups of \textsf{GAP}~\cite{GAP4-11-1}) and to $C_2\times C_2\times C_4\times C_4$ (\texttt{SmallGroup(64, 192)}). The isomorphism classes of left braces with additive group isomorphic to $C_4\times C_4\times C_4$ have been obtained recently by the authors  in \cite{BallesterEstebanPerezC-64-55}. In this paper we analyse the remaining case and we conclude the classification of left braces of order~$64$.

\begin{theorem}\label{th-A}
  The number of isomorphism classes of left braces of order~$64$ with additive group  isomorphic to $C_2\times C_2\times C_4\times C_4$ is  $\numberbraces$.
\end{theorem}

As we have commented in \cite{BallesterEstebanPerezC-64-55}, the determination of the isomorphism classes of left braces with additive group isomorphic to a given abelian group~$G$ reduces to the determination of all conjugacy classes of regular subgroups of the holomorph of~$G$. A detailed study of the holomorph of a group is presented in~\cite{BallesterEstebanPerezC-64-55}.

The techniques that we use here to obtain the conjugacy classes of regular subgroups of the holomorph $C_2\times C_2\times C_4\times C_4$ differ from the ones used in~\cite{BallesterEstebanPerezC-64-55} to get the conjugacy classes of regular subgroups of the holomorph of $C_4\times C_4\times C_4$. The main reasons for using a different algorithm are that in our computation we obtain a large number of intermediate subgroups, larger than the ones in~\cite{BallesterEstebanPerezC-64-55} (see Table~\ref{tab-layer} below), and that we are finding that the read-write operations on files are requiring a considerable amount of time. We have considered that a different approach that combines saving some information to files with minimising the read-write operations might be more convenient. We hope that the ideas in this note could be useful to determine the isomorphism classes of skew left braces with other additive groups.

As a consequence of Theorem~\ref{th-A} together with the results of \cite{BardakovNeshchadimYadav20} and  \cite{BallesterEstebanPerezC-64-55}, we can now give the number of isomorphism classes of left braces of order~$64$.
\begin{corollary}
 The number of isomorphism classes of left braces of order~$64$ is $\totalbraces$.
\end{corollary}
\begin{table}
  \centering
  \caption{Left braces of order~$64$ by additive group}
  \label{tab-braces-64-BNY20}
  \begin{tabular}{llr}\hline
    Group id&Structure&Number\\\hline
    $(64, 1)$&$C_{64}$&10\\
    $(64, 2)$&$C_8\times C_8$&$11\,354$\\
    $(64, 26)$&$C_{4}\times C_{16}$& $2\,724$\\
    $(64, 50)$&$C_{2}\times C_{32}$&$142$\\
    $(64,55)$&$C_4\times C_4\times C_4$&$\numberbracesold$\\
    $(64,83)$&$C_2\times C_4\times C_8$&$743\,410$\\
    $(64, 183)$&$C_2\times C_2\times C_{16}$&$3\,124$\\
    $(64, 192)$&$C_2\times C_2\times C_4\times C_4$&$\numberbraces$\\
    $(64, 246)$&$C_2\times C_2\times C_2\times C_8$&$253\,350$\\
    $(64, 260)$&$C_2\times C_2\times C_2\times C_2\times C_4$&$2\,189\,661$\\
    $(64, 267)$&$C_2\times C_2\times C_2\times C_2\times C_2\times C_2$&$58\,558$\\\hline
    Total & & $\totalbraces$\\\hline
  \end{tabular}
\end{table}

Table~\ref{tab-braces-64-BNY20} records the numbers of braces of order~$64$ by additive group given in~\cite[Table~3]{BardakovNeshchadimYadav20}, in~\cite{BallesterEstebanPerezC-64-55}, and in Theorem~\ref{th-A} by additive group.
Table~\ref{tab-mult} summarises the isomorphism classes of braces with additive group $C_2\times C_2\times C_4\times C_4$ by the isomorphism class of the multiplicative group. For instance, the first two entries of the first row of the table, namely $4$, $185$, mean that there are $185$ isomorphism classes of braces with additive group $C_2\times C_2\times C_4\times C_4$ and multiplicative group isomorphic to \texttt{SmallGroup(64, 4)}.
\begin{table}[htbp]
  \caption{Number of isomorphism classes of left braces with additive group $C_2\times C_2\times C_4\times C_4$ by multiplicative group}\label{tab-mult}\medskip
  \footnotesize
  \centering
    \begin{tabular}{rrrrrrrrrrrr}\hline
    id&\#&id&\#&id&\#&id&\#&id&\#&id&\#\\\hline
    4 & 185&     69 & 163\,216&  103 & 48&      146 & 262&     200 & 19\,003&  234 & 504\,164\\ 
    5 & 89&      70 & 80\,652&   104 & 39&      147 & 64&      201 & 81\,305&  235 & 233\,582\\ 
    7 & 4&       71 & 44\,046&   105 & 44&      148 & 198&     202 & 36\,068&  236 & 126\,254\\ 
    8 & 152&     72 & 40\,622&   106 & 40&      149 & 262&     203 & 96\,180&  237 & 233\,576\\ 
    9 & 218&     73 & 31\,133&   107 & 40&      150 & 64&      204 & 83\,536&  238 & 71\,696\\  
    10 & 2&      74 & 27\,623&   108 & 46&      151 & 198&     205 & 90\,204&  239 & 12\,016\\  
    14 & 2&      75 & 86\,631&   109 & 146&     156 & 56&      206 & 264\,425& 240 & 126\,305\\ 
    17 & 47&     76 & 27\,293&   112 & 70&      158 & 56&      207 & 44\,746&  241 & 181\,022\\ 
    18 & 236&    77 & 55\,112&   113 & 68&      159 & 2&       208 & 38\,531&  242 & 34\,059\\  
    20 & 399&    78 & 80\,755&   114 & 60&      160 & 116&     209 & 54\,734&  243 & 252\,024\\ 
    21 & 80&     79 & 79\,679&   115 & 132&     161 & 92&      210 & 485\,604& 244 & 234\,017\\ 
    23 & 2\,709&   80 & 27\,603&   116 & 264&     162 & 90&      211 & 8\,857&   245 & 14\,338\\  
    24 & 78&     81 & 77\,603&   117 & 136&     163 & 130&     212 & 19\,738&  246 & 40\\     
    25 & 232&    82 & 7\,208&    119 & 122&     164 & 170&     213 & 131\,810& 247 & 169\\    
    32 & 852&    83 & 25&      120 & 122&     165 & 168&     214 & 112\,426& 248 & 108\\    
    33 & 730&    84 & 23&      121 & 120&     166 & 254&     215 & 143\,656& 249 & 122\\    
    34 & 372&    85 & 39&      122 & 120&     168 & 56&      216 & 142\,162& 250 & 30\\     
    35 & 410&    86 & 44&      128 & 129&     169 & 58&      217 & 122\,955& 251 & 118\\    
    36 & 96&     87 & 319&     129 & 275&     170 & 56&      218 & 131\,619& 252 & 88\\     
    37 & 102&    88 & 344&     130 & 334&     172 & 56&      219 & 522\,728& 253 & 328\\    
    55 & 1\,054&   89 & 274&     131 & 257&     173 & 32&      220 & 485\,625& 254 & 386\\    
    56 & 8\,754&   90 & 6\,647&    132 & 403&     175 & 32&      221 & 261\,298& 255 & 556\\    
    57 & 6\,624&   91 & 2\,785&    133 & 554&     176 & 64&      222 & 224\,308& 256 & 648\\    
    58 & 42\,614&  92 & 439&     134 & 560&     177 & 32&      223 & 485\,594& 257 & 156\\    
    59 & 40\,408&  93 & 161&     135 & 232&     178 & 96&      224 & 40\,841&  258 & 444\\    
    60 & 8\,897&   94 & 331&     136 & 700&     180 & 2&       225 & 112\,213& 259 & 288\\    
    61 & 42\,552&  95 & 144&     137 & 352&     192 & 3\,972&    226 & 147\,895& 260 & 978\\    
    62 & 21\,370&  96 & 232&     138 & 4\,656&    193 & 13\,622&   227 & 546\,020& 261 & 5\,564\\   
    63 & 19\,782&  97 & 290&     139 & 4\,329&    194 & 11\,511&   228 & 253\,872& 262 & 1\,316\\   
    64 & 12\,778&  98 & 336&     141 & 64&      195 & 42\,161&   229 & 168\,878& 263 & 16\,329\\  
    65 & 7\,034&   99 & 176&     142 & 120&     196 & 89\,797&   230 & 78\,116&  264 & 17\,318\\  
    66 & 85\,060&  100 & 270&    143 & 184&     197 & 25\,927&   231 & 90\,758&  265 & 6\,987\\   
    67 & 89\,328&  101 & 742&    144 & 120&     198 & 81\,586&   232 & 504\,784& 266 & 10\,268\\  
    68 & 157\,112& 102 & 713&    145 & 250&     199 & 67\,162&   233 & 467\,180& 267 & 46\\\hline         
  \end{tabular}
\end{table}

\section{The computations}


As we have mentioned in the introduction, the determination of the isomorphism classes of left braces with additive group $G=C_2\times C_2\times C_4\times C_4$ is reduced to the computation of all conjugacy classes of regular subgroups of~$G$. The advantage in this case with respect to the case  $C_4\times C_4\times C_4$ analysed  in \cite{BallesterEstebanPerezC-64-55} is that the holomorph of $G$ is soluble, whereas the holomorph of $C_4\times C_4\times C_4$ was insoluble. This makes suitable the application of an adaptation of the algorithm of Hulpke \cite{Hulpke99} to determine the conjugacy classes of subgroups of a soluble group $S$, implemented in the function \texttt{SubgroupsSolvableGroup} of \textsf{GAP} \cite{GAP4-11-1}. This algorithm is based on a normal series $1=N_r\le N_{r-1}\dots\le N_2\le N_1\le N_0=S$ of a soluble group $S$ whose factors are elementary abelian. The method consists on the computation by induction of the subgroups of $S/N_{i+1}$, that we will say to belong to the layer $i+1$, from the subgroups of $S/N_i$ previously computed, that we will say to belong to the layer $i$, in the following way. Suppose that $U/N_{i+1}$ is a subgroup of $S/N_{i+1}$. Then $UN_i/N_i$ is a subgroup of $S/N_i$ and so we have one of the following possibilities for $UN_i/N_i$:
\begin{enumerate}
\item $N_i< U$ and so $UN_i/N_{i+1}=U/N_{i+1}$ can be taken to be the preimage of $U/N_i$ under the natural epimorphism from $S/N_{i+1}$ to $S/N_{i}$.
\item $U\le N_i$ and so $U/N_{i+1}$ is a subgroup of the elementary abelian group $N_i/N_{i+1}$.
\item $U\not\le N_i$ and $N_i\not \le U$. In this case, $UN_i/N_i$, the image of $U/N_{i+1}$ under the natural epimorphism from $S/N_{i+1}$ to $S/N_i$, is a non-trivial subgroup of $S/N_i$. Furthermore, $B/N_{i+1}:=(U/N_{i+1})\cap (N_i/N_{i+1})=(U\cap N_i)/N_{i+1}$ is a proper subgroup of the elementary abelian group $N_i/N_{i+1}$ and $(U/N_{i+1})/(B/N_{i+1})$ complements $(N_i/N_{i+1})/(B/N_{i+1})$ in $(UN_i/N_{i+1})/(B/N_{i+1})$.
\end{enumerate}

  

The implementation of \texttt{SubgroupsSolvableGroup} in \textsf{GAP} allows adding restrictions like \texttt{ExactSizeConsiderFunction} to avoid the computation of subgroups that do not lead to subgroups of the specified order. This is useful since regular subgroups of $\Hol(G)$ have the same order as $G$. We also note that regular subgroups of $\Hol(G)$ must have a surjective ``projection'' onto $G$ and so we can add the restriction of \cite[Proposition~2.3]{BallesterEstebanPerezC-64-55} to discard all subgroups leading only to non-regular subgroups.

The \textsf{GAP} implementation of \texttt{SubgroupsSolvableGroup} seems to prefer series with elementary abelian factors with large factor groups for the computation of conjugacy classes of subgroups, especially the derived series whenever this series has elementary abelian factors. However, in this case we have obtained that the number of subgroups of the elementary abelian factors $N_i/N_{i+1}$ can be huge and even the computation of the orbits of these subgroups with respect to the action of the holomorph of $G$ can be very tedious, as well as the determination of conjugacy classes of complements. Nevertheless, the algorithm also allows to feed the calculation with a normal series with elementary abelian factors. Hence we have decided to work with a series of length~$11$ in which all factors below the derived subgroup are chief factors. This makes the computation of conjugacy classes of subgroups of each factor group fast as well as the determination of complements in the groups obtained in the previous layers. However, we face another problem: the numbers of subgroups obtained at some layers are growing very fast, as we can see on the second column of Table~\ref{tab-layer}.
\begin{table}
  \centering
  \caption{Numbers of conjugacy classes of subgroups obtained and jobs used at each layer}\label{tab-layer}
  \begin{tabular}{rrr}\hline
    Layer&\# classes&\# jobs\\\hline
    2&5&1\\
    3&22&2\\
    4&127&3\\
    5&17\,025&4\\
    6&5\,644\,351&5\\
    7&5\,796\,845&562\\
    8&8\,712\,963&562\\
    9&90\,704\,762&812\\
    10&1\,593\,095\,679&1\,505\\
    11&10\,326\,821&1\,505\\\hline
  \end{tabular}
\end{table}

As we have mentioned in \cite{BallesterEstebanPerezC-64-55}, in the \textsf{GAP} implementation of the function \texttt{SubgroupsSolvableGroup}, all conjugacy classes of subgroups of $G/N_i$ (layer~$i$) and all computed conjugacy classes of $G/N_{i+1}$ (layer~$i+1$) are stored at each layer. However, we only need one group of the layer~$i$ is at each step and we will not use the groups obtained for the layer $i+1$ until advancing to the next layer. Furthermore, when there are a large number of subgroups, they can use a huge amount of memory, that might eventually exhaust the physical memory. In our case, it seems that the number of subgroups in layer~$10$, even if we find a cheap representation of a subgroup in terms of memory, would lead to a huge memory usage, probably beyond what most personal computers with up to 32 or 64~Gb of memory can handle. Our solution in this case has been to obtain the subgroups of the elementary abelian group $N_{i-1}/N_i$ for each layer $i$ and to decompose the list of subgroups of $G/N_{i-1}$ obtained  at the layer $i-1$ into several sublists and to compute the subgroups corresponding to the layer $i$ for each of these sublists. This operation can be clearly done in parallel once the subgroups of the elementary abelian group $N_{i-1}/N_i$ are known. The groups obtained at each layer for each sublist can be saved on the hard disk and the obtained list or a sublist of this list can be used as an input for the computation in the next layer. To do so, we store a generating set of each subgroup in such a list as an integer by using the same ideas of the function \texttt{CodePcGroup} of \textsf{GAP} for a generating set instead of the power-commutator relations.

Once we had noted during the analysis of layer~$9$, with $90\,704\,762$ subgroups, that the sequential implementation used in \cite{BallesterEstebanPerezC-64-55} could need a lot of computing time and disk access time on a computer with an Intel processor i7-11700 and 32~Gb of RAM, we have implemented these computations in the supercomputer \emph{Llu\'\i s Vives} of the Universitat de Val\`encia \cite{LluisVivesv2}. We thank the Computer Service for granting us an immediate access to this facility and for helping us by installing \textsf{GAP} and solving all our queries. Our implementation includes a modification of Hulpke's algorithm that saves the subgroups of $N_{i-1}/N_i$ that can be obtained at a given layer and that computes only the subgroups that appear at a given layer $i$, corresponding to $G/N_i$, from a list of subgroups of $G/N_{i-1}$ obtained in the calculation of the previous layer. Our code includes a restriction to discard the intermediate subgroups obtained at each layer that cannot contain regular subgroups because the corresponding ``projections'' on the additive group are not surjective. We have also included the ``consider'' function \texttt{ExactSizeConsiderFunction(64)} and the series with elementary abelian factors we have described before. The results of the computation have been obtained in less than two days with some peaks of more than 550 simultaneous jobs. The last column of Table~\ref{tab-layer} shows the number of jobs we have used for each layer, which coincides with the number of sublists in which we have decomposed the list of groups of the previous layer to get the groups in the current layer.

The algorithm was tested before with the search of the conjugacy classes of regular subgroups of the holomorphs of $C_4\times C_{16}$ (\texttt{SmallGroup(64, 26)}), $C_2\times C_4\times C_8$ (\texttt{SmallGroup(64, 83)}), and $C_2\times C_2\times C_{16}$ (\texttt{SmallGroup(64, 183)}). The numbers of conjugacy classes obtained by our algorithm coincided with the ones obtained in \cite{BardakovNeshchadimYadav20} and that appear in Table~\ref{tab-braces-64-BNY20}.

\section*{Acknowledgements}

We thank the Computer Service of the Universitat de Val\`encia for allowing us to use the scientific supercomputer \emph{Llu\'\i s Vives} and for their kind support in setting up the system for the computation and solving our questions.

\section*{Declarations}
\subsection*{Funding}
These results are part of the R+D+i project supported by the Grant 
PGC2018-095140-B-I00, funded by MCIN/AEI/10.13039/501100011033 and by ``ERDF A way of making Europe''.

\subsection*{Data availability}
The list of left braces with additive group $C_2\times C_2\times C_4\times C_4$ is available for use in \textsf{GAP} \cite{GAP4-11-1} on \url{https://github.com/RamonEstebanRomero/braces64} \cite{BallesterEstebanPerezC22-braces64-GitHub}, together with the rest of the left braces of order~$64$. The storage of these left braces is based on the ideas of \cite{BallesterEsteban22} of interpreting them as trifactorised groups. We also include in this repository also some \textsf{GAP} functions that allow us to work with these left braces with the help of the \textsf{YangBaxter} package \cite{VendraminKonovalov22-YangBaxter-0.10.1} for \textsf{GAP}. The \textsf{GAP} code used for the computation is available from the authors on request.

\subsection*{Conflict of interest}

On behalf of all authors, the corresponding author states there is no conflict of interest.

\bibliographystyle{plain}
\bibliography{bibgroup}
\end{document}